\documentclass[10pt,twoside,reqno]{amsart}
\usepackage{amssymb}
\usepackage{multirow}
\calclayout

\numberwithin{equation}{section}
\makeatletter

\renewcommand{\@secnumfont}{\bfseries}

\renewcommand{\section}{\@startsection{section}{1}%
  {0mm}{.7\linespacing\@plus\linespacing}{.5\linespacing}%\z@={0mm}
  {\normalfont\bfseries\centering}}

\newcommand{\bibsection}{\@startsection{section}{1}%
  {0mm}{.7\linespacing\@plus\linespacing}{.5\linespacing}%\z@={0mm}
  {\normalfont\scshape\centering}}

\renewcommand{\@biblabel}[1]{#1.}

\newtheorem{thm}{\bf Theorem}[section]

\newtheorem{cor}[thm]{\bf Corollary}

\newtheorem{rmk}{\bf Remark}

\begin{document}

\vspace{1.3cm}

\title {Degenerate Bernstein polynomials}

\author{Taekyun Kim}
\address{Department of Mathematics, Kwangwoon University, Seoul 139-701, Republic of Korea}
\email{tkkim@kw.ac.kr}

\author{Dae San Kim}
\address{Department of Mathematics, Sogang University, Seoul 121-742, Republic of Korea}
\email{dskim@sogang.ac.kr}


\begin{abstract} 
Here we consider the degenerate Bernstein polynomials as a degenerate version of Bernstein polynomials, which are motivated by Simsek's recent work 'Generating functions for unification of the multidimensional Bernstein polynomials and their applications'([15,16]) and Carlitz's degenerate Bernoulli polynomials. We derived thier generating function, symmetric identities, recurrence relations, and some connections with generalized falling factorial polynomials, higher-order degenerate Bernoulli polynomials and degenerate Stirling numbers of the second kind.
\end{abstract}

\maketitle

\section{Introduction}

For $\lambda \in \mathbb{R}$, the degenerate Bernoulli polynomials of order $k$ are defined by L. Cartliz as
\begin{equation}\begin{split}\label{01}
\left( \frac{t}{(1+\lambda t)^{\frac{1}{\lambda }}-1} \right)^k (1+\lambda t)^{\frac{x}{\lambda }} = \sum_{n=0}^\infty  \beta_{n,\lambda }^{(k)}      \frac{t^n}{n!},\quad (\text{see}\,\, [4,5]).
\end{split}\end{equation}

Note that $\lim_{\lambda \rightarrow 0} \beta_{n,\lambda }^{(k)} (x) = B_n^{(k)}(x)$ are the ordinary Bernoulli polynomials of order $k$ given by 
\begin{equation*}\begin{split}
\left( \frac{t}{e^t-1}\right)^k e^{xt} = \sum_{n=0}^\infty    B_n^{(k)}(x)    \frac{t^n}{n!},\quad (\text{see}\,\, [1,12,13]).
\end{split}\end{equation*}

It is known that the Stirling number of the second kind is defined as
\begin{equation}\begin{split}\label{02}
x^n = \sum_{l=0}^n S_2(n,l) (x)_l,\quad (\text{see}\,\, [2,4,8,10]),
\end{split}\end{equation}
where $(x)_l = x(x-1) \cdots (x-l+1)$, $(l \geq 1),$ $(x)_0=1$.

For $\lambda \in \mathbb{R}$, the $(x)_{n,\lambda}$ is defined as
\begin{equation}\begin{split}\label{03}
(x)_{0,\lambda }=1,\,\,(x)_{n,\lambda }=x(x-\lambda )(x-2\lambda )\cdots(x-(n-1)\lambda ),\,\,(n \geq 1)
\end{split}\end{equation}

In [8,9,10], ${x \choose n}_\lambda$ is defined as
\begin{equation}\begin{split}\label{04}
{x \choose n}_\lambda  = \frac{(x)_{n,\lambda }}{n!} = \frac{x(x-\lambda )\cdots(x-(n-1)\lambda )}{n!},\,\,(n \geq 1),\,\, {x \choose 0}_\lambda  =1.
\end{split}\end{equation}

Thus, by \eqref{04}, we get
\begin{equation}\begin{split}\label{05}
(1+\lambda t)^{\frac{x}{\lambda }} = \sum_{n=0}^\infty {x \choose n}_\lambda  t^n,\quad (\text{see}\,\, [7]).
\end{split}\end{equation}

From \eqref{05}, we note that
\begin{equation}\begin{split}\label{06}
\sum_{m=0}^n {y \choose m}_\lambda  {x \choose n-m}_\lambda  = {x+y \choose n}_\lambda ,\,\,(n \geq 0).
\end{split}\end{equation}

The degenerate Stirling numbers of the second kind are defined by
\begin{equation}\begin{split}\label{07}
\frac{1}{k!} \big( (1+\lambda t)^{\frac{1}{\lambda }}-1\big)^k = \sum_{n=k}^\infty S_{2,\lambda }(n,k) \frac{t^n}{n!},\,\,(k \geq 0),\quad (\text{see}\,\, [7,8]).
\end{split}\end{equation}

By \eqref{07}, we easily get
\begin{equation*}\begin{split}
\lim_{\lambda  \rightarrow 0} S_{2,\lambda }(n,k) = S_2(n,k),\,\,(n \geq k \geq 0),\quad (\text{see}\,\, [8,10]).
\end{split}\end{equation*}

In this paper, we use the following notation.
\begin{equation}\begin{split}\label{08}
(x \oplus_\lambda  y)^n = \sum_{k=0}^n {n \choose k} (x)_{k,\lambda } (y)_{n-k,\lambda },\,\,(n \geq 0).
\end{split}\end{equation}

The Bernstein polynomials of degree $n$ is defined by
\begin{equation}\begin{split}\label{09}
B_{k,n}(x) = {n \choose k} x^k (1-x)^{n-k},\,\,(n \geq k \geq 0),\quad (\text{see}\,\, [6,11,17]).
\end{split}\end{equation}

Let $C[0,1]$ be the space of continuous functions on $[0,1]$. The Bernstein operator of order $n$ for $f$ is given by
\begin{equation}\begin{split}\label{10}
\mathbb{B}_n (f|x) = \sum_{k=0}^n f (\tfrac{k}{n} ) {n \choose k} x^k (1-x)^{n-k} = \sum_{k=0}^n  f (\tfrac{k}{n} ) B_{k,n}(x),
\end{split}\end{equation}
where $n \in \mathbb{N} \cup \{0\}$ and $f \in C[0,1]$, (see [3,6,14]). 

A Bernoulli trial involves performing a random experiment and noting whether a particular event A occurs.
The outcome of Bernoulli trial is said to be "success" if A occurs and a "failure" otherwise. The probability $P_n(k)$ of $k$ successes in $n$ independent Bernoulli trials is given by the binomial probability law:
\begin{equation*}\begin{split}
P_n (k) = {n \choose k} p^k (1-p)^{n-k},\,\,\text{for}\,\,k=0,1,2,\cdots,
\end{split}\end{equation*}
From the definition of Bernstein polynomials we note that Bernstein basis is probability mass of binomial distribution with parameter $(n, x=p)$. 

Let us assume that the probability of success in an experiment is $p$. We wondered if we can say the probability of success in the nineth trial is still $p$ after failing eight times in a ten trial experiment.  Because there's a psychological burden to be successful. 

It seems plausible that the probability is less than $p$. This speculation motivated the study of the degenerate Bernstein polynomials associated with the probability distribution.

In this paper, we consider the degenerate Bernstein polynomials as a degenerate version of Bernstein polynomials. We derive thier generating function, symmetric identities, recurrence relations, and some connections with generalized falling factorial polynomials, higher-order degenerate Bernoulli polynomials and degenerate Stirling numbers of the second kind.

\section{Degenerate Bernstein polynomials}

For $\lambda \in \mathbb{R}$ and $k,n \in \mathbb{N}\cup \{0\}$, with $k \leq n$, we define the degenerate Bernstein polynomials of degree $n$ which are given by
\begin{equation}\begin{split}\label{11}
B_{k,n}(x|\lambda ) = {n \choose k} (x)_{k,\lambda }(1-x)_{n-k,\lambda },\,\,(x \in [0,1]).
\end{split}\end{equation}

Note that $\lim_{\lambda \rightarrow 0} B_{k,n}(x|\lambda ) = B_{k,n}(x)$, $(0 \leq k \leq n)$.
From \eqref{11}, we derive the generating function of $B_{k,n}(x|\lambda )$, which are given by
\begin{equation}\begin{split}\label{12}
\sum_{n=k}^\infty B_{k,n}(x|\lambda ) \frac{t^n}{n!} &= \sum_{n=k}^\infty {n \choose k} (x)_{k,\lambda }(1-x)_{n-k,\lambda } \frac{t^n}{n!}\\
&= \frac{(x)_{k,\lambda}}{k!} \sum_{n=k}^\infty \frac{1}{(n-k)!} (1-x)_{n-k,\lambda }t^n\\
&= \frac{(x)_{k,\lambda }}{k!} \sum_{n=0}^\infty \frac{(1-x)_{n,\lambda }}{n!} t^{n+k} \\
&= \frac{(x)_{k,\lambda }}{k!} t^k \sum_{n=0}^\infty {1-x \choose n}_\lambda  t^n \\
&= \frac{(x)_{k,\lambda }}{k!}  t^k (1+\lambda t)^{\frac{1-x}{\lambda }}.
\end{split}\end{equation}

Therefore, by \eqref{12}, we obtain the following theorem.
\begin{thm}
For $x \in [0,1]$ and $k=0,1,2,\cdots,$ we have
\begin{equation*}\begin{split}
\frac{1}{k!} (x)_{k,\lambda } t^k (1+\lambda t)^{\frac{1-x}{\lambda }} = \sum_{n=k}^\infty B_{k,n}(x|\lambda ) \frac{t^n}{n!}.
\end{split}\end{equation*}
\end{thm}

From \eqref{11}, we note that
\begin{equation}\begin{split}\label{13}
B_{k,n}(x|\lambda ) = {n \choose k} (x)_{k,\lambda } (1-x)_{n-k,\lambda } = {n \choose n-k} (x)_{k,\lambda }(1-x)_{n-k,\lambda }.
\end{split}\end{equation}

By replacing $x$ by $1-x$, we get
\begin{equation}\begin{split}\label{14}
B_{k,n}(1-x|\lambda ) = {n \choose n-k} (1-x)_{k,\lambda }(x)_{n-k,\lambda } = B_{n-k,n}(x|\lambda ),
\end{split}\end{equation}
where $n,k \in \mathbb{N} \cup \{0\}$, with $0 \leq k \leq n$.

Therefore, by \eqref{14}, we obtain the following theorem.
\begin{thm} (\textnormal{Symmetric identities})
For $n,k \in \mathbb{N} \cup \{0\}$, with $k \leq n$, and $x \in [0,1]$, we have
\begin{equation*}\begin{split}
B_{n-k,n}(x|\lambda )=B_{k,n}(1-x|\lambda )  .
\end{split}\end{equation*}
\end{thm}

Now, we observe that
\begin{equation}\begin{split}\label{15}
&\frac{n-k}{n} B_{k,n}(x|\lambda ) + \frac{k+1}{n} B_{k+1,n}(x|\lambda )\\
&= \frac{n-k}{n} {n \choose k} (x)_{k,\lambda } (1-x)_{n-k,\lambda } + \frac{k+1}{n} {n \choose k+1} (x)_{k+1,\lambda }(1-x)_{n-k-1,\lambda }\\
&= \frac{(n-1)!}{k! (n-k-1)!} (x)_{k,\lambda }(1-x)_{n-k,\lambda } + \frac{(n-1)!}{k! (n-k-1)!} (x)_{k+1,\lambda } (1-x)_{n-k-1,\lambda }\\
&=(1-x-(n-k-1)\lambda ) B_{k,n-1}(x|\lambda ) + (x-k\lambda )B_{k,n-1}(x|\lambda )\\
&= (1+ \lambda (1-n))B_{k,n-1}(x|\lambda).
\end{split}\end{equation}

Therefore, by \eqref{15}, we obtain the following theorem.
\begin{thm}
For $k \in \mathbb{N} \cup \{0\}$, $n \in \mathbb{N}$, with $k \leq n-1$, and $x \in [0,1]$, we have
\begin{equation}\begin{split}\label{16}
(n-k)B_{k,n}(x|\lambda ) + (k+1)B_{k+1,n}(x|\lambda )=(1+ \lambda (1-n))B_{k,n-1}(x|\lambda).
\end{split}\end{equation}
\end{thm}

From \eqref{11}, we have
\begin{equation}\begin{split}\label{17}
&\left( \frac{n-k+1}{k} \right) \left( \frac{n-(k-1)\lambda }{1-x-(n-k)\lambda } \right) B_{k-1,n}(x|\lambda ) \\
&=\left( \frac{n-k+1}{k} \right) \left( \frac{n-(k-1)\lambda }{1-x-(n-k)\lambda } \right)
{n \choose k-1} (x)_{k-1,\lambda }(1-x)_{n-k+1,\lambda }\\
&= \frac{n!}{k! (n-k)!} (x)_{k-1,\lambda }(1-x)_{n-k,\lambda } = B_{k,n}(x|\lambda ).
\end{split}\end{equation}

Therefore, by \eqref{17}, we obtain the following theorem.
\begin{thm} For $n,k \in \mathbb{N}$, with $k \leq n$, we have
\begin{equation*}\begin{split}
\left( \frac{n-k+1}{k} \right) \left( \frac{n-(k-1)\lambda }{1-x-(n-k)\lambda } \right) B_{k-1,n}(x|\lambda ) = B_{k,n}(x|\lambda ).
\end{split}\end{equation*}
\end{thm}

For $0 \leq k \leq n$, we get
\begin{equation}\begin{split}\label{18}
&(1-x-(n-k-1)\lambda ) B_{k,n-1}(x|\lambda ) + (x-(k-1)\lambda )B_{k-1,n-1}(x|\lambda ) \\
&=(1-x-(n-k-1)\lambda ) {n-1 \choose k} (x)_{k,\lambda }(1-x)_{n-1-k,\lambda } \\
&\quad + (x-(k-1)\lambda ){n-1 \choose k-1}(x)_{k-1,\lambda }(1-x)_{n-k,\lambda }\\
&={n-1 \choose k} (x)_{k,\lambda }(1-x)_{n-k,\lambda }+ {n-1 \choose k-1} (x)_{k,\lambda } (1-x)_{n-k,\lambda }\\
&=\left( {n-1 \choose k}+{n-1 \choose k-1} \right) (x)_{k,\lambda } (1-x)_{n-k,\lambda } = {n \choose k} (x)_{k,\lambda } (1-x)_{n-k,\lambda }.
\end{split}\end{equation}

Therefore, by \eqref{18}, we obtain the following theorem.
\begin{thm}(\textnormal{Recurrence formula}).
For $k,n \in \mathbb{N}$, with $k \leq n-1$, $x \in [0,1]$, we have
\begin{equation*}\begin{split}
(1-x-(n-k-1)\lambda ) B_{k,n-1}(x|\lambda ) + (x-(k-1)\lambda )B_{k-1,n-1}(x|\lambda )=B_{k,n}(x|\lambda ).
\end{split}\end{equation*}
\end{thm}

\begin{rmk}For $n\in \mathbb{N}$, we have
\begin{equation*}\begin{split}
&\sum_{k=0}^n \frac{k}{n} B_{k,n}(x|\lambda ) = \sum_{k=0}^n \frac{k}{n} {n \choose k} (x)_{k,\lambda } (1-x)_{n-k,\lambda }\\
&= \sum_{k=1}^n {n-1 \choose k-1} (x)_{k,\lambda } (1-x)_{n-k,\lambda } = \sum_{k=0}^{n-1} {n-1 \choose k} (x)_{k+1,\lambda }(1-x)_{n-1-k,\lambda }\\
&=(x-k\lambda ) \sum_{k=0}^{n-1} {n-1 \choose k} (x)_{k,\lambda } (1-x)_{n-1-k,\lambda } = (x-k\lambda ) (x \oplus_\lambda  (1-x) )^{n-1}.
\end{split}\end{equation*}
\end{rmk}

Now, we observe that
\begin{equation}\begin{split}\label{19}
\sum_{k=2}^n \frac{{k \choose 2}}{{n \choose 2}} B_{k,n}(x|\lambda ) &=
\sum_{k=2}^n \frac{k(k-1)}{n(n-1)} {n \choose k} (x)_{k,\lambda }(1-x)_{n-k,\lambda }\\
&= \sum_{k=2}^n \frac{k(k-1)}{n(n-1)} {n \choose k} (x)_{k,\lambda }(1-x)_{n-k,\lambda }\\
&=\sum_{k=2}^n {n-2 \choose k-2} (x)_{k,\lambda } (1-x)_{n-k,\lambda }\\
&= \sum_{k=0}^{n-2} {n-2 \choose k} (x)_{k+2,\lambda }(1-x)_{n-2-k,\lambda }\\
&= (x-k\lambda )(x-(k+1)\lambda ) \sum_{k=0}^{n-2} {n-2 \choose k} (x)_{k,\lambda } (1-x)_{n-2-k,\lambda }.
\end{split}\end{equation}

Similarly, we have
\begin{equation}\begin{split}\label{20}
\sum_{k=i}^n \frac{{k \choose i}}{{n \choose i}} B_{k,n}(x|\lambda ) &= (x-k\lambda )_{i,\lambda } \sum_{k=0}^{n-i} {n-i \choose k} (x)_{k,\lambda } (1-x)_{n-i-k,\lambda }\\
&= (x-k\lambda )_{i,\lambda } (x \oplus_\lambda  (1-x))^{n-i}.
\end{split}\end{equation}

From \eqref{20}, we note that
\begin{equation}\begin{split}\label{21}
(x-k\lambda)_{i,\lambda } = \frac{1}{(x \oplus_\lambda  (1-x))^{n-i}} \sum_{k=i}^n \frac{{k \choose i}}{{n \choose i}} B_{k,n}(x|\lambda ),
\end{split}\end{equation}
where $n,i \in \mathbb{N}$, with $ i \leq n$, and $x \in [0,1]$.

Therefore, by \eqref{21}, we obtain the following theorem.
\begin{thm}
For $n,i \in \mathbb{N}$, with $ i \leq n$, and $x \in [0,1]$, we have
\begin{equation*}\begin{split}
(x-k\lambda)_{i,\lambda } = \frac{1}{(x \oplus_\lambda  (1-x))^{n-i}} \sum_{k=i}^n \frac{{k \choose i}}{{n \choose i}} B_{k,n}(x|\lambda ).
\end{split}\end{equation*}
\end{thm}

From Theorem 2.1, we note that
\begin{equation}\begin{split}\label{22}
&\frac{t^k}{k!} (x)_{k,\lambda }(1+\lambda t)^{\frac{1-x}{\lambda }} = \frac{(x)_{k,\lambda }}{k!} \Big( (1+\lambda t)^{\frac{1}{\lambda }}-1 \Big)^k \left(\frac{t}{(1+\lambda t)^{\frac{1}{\lambda }}-1} \right)^k (1+\lambda t)^{\frac{1-x}{\lambda }}\\
&= (x)_{k,\lambda } \left( \sum_{m=k}^\infty S_{2,\lambda }(m,k) \frac{t^m}{m!} \right) \left( \sum_{l=0}^\infty \beta_{l,\lambda }^{(k)} (1-x) \frac{t^l}{l!} \right)\\
&= (x)_{k,\lambda } \sum_{n=k}^\infty \left( \sum_{m=k}^n {n \choose m} S_{2,\lambda }(m,k) \beta_{n-m,\lambda }^{(k)} (1-x) \right) \frac{t^n}{n!}.
\end{split}\end{equation}

On the other hand,
\begin{equation}\begin{split}\label{23}
\frac{(x)_{k,\lambda }}{k!} t^k (1+\lambda t)^{\frac{1-x}{\lambda }} = \sum_{n=k}^\infty B_{k,n}(x|\lambda ) \frac{t^n}{n!}.
\end{split}\end{equation}

Therefore, by \eqref{22} and \eqref{23}, we obtain the following theorem.

\begin{thm}
For $n,k \in \mathbb{N}\cup \{0\}$ with $n \geq k$, we have
\begin{equation*}\begin{split}
 B_{k,n}(x|\lambda ) = (x)_{k,\lambda } \sum_{m=k}^n {n \choose m} S_{2,\lambda }(m,k) \beta_{n-m,\lambda }^{(k)} (1-x).
 \end{split}\end{equation*} 
\end{thm}

Let $\Delta$ be the shift difference operator with $\Delta f(x) = f(x+1)-f(x)$. Then we easily get
\begin{equation}\begin{split}\label{24}
\Delta^n f(0) = \sum_{k=0}^n {n \choose k} (-1)^{n-k} f(k), \,\,(n\in \mathbb{N}\cup \{0\}).
\end{split}\end{equation}

Let us take $f(x)=(x)_{m,\lambda }$, $(m \geq 0)$. Then, by \eqref{24}, we get
\begin{equation}\begin{split}\label{25}
\Delta^n (0)_{m,\lambda } = \sum_{k=0}^n {n \choose k} (-1)^{n-k} (k)_{m,\lambda }.
\end{split}\end{equation}

From \eqref{07}, we note that
\begin{equation}\begin{split}\label{26}
\sum_{n=k}^\infty S_{2,\lambda }(n,k) \frac{t^n}{n!}&=
\frac{1}{k!} \Big( (1+\lambda t)^{\frac{1}{\lambda }}-1 \Big)^k = \frac{1}{k!} \sum_{l=0}^k {k \choose l} (-1)^{k-l} (1+\lambda t)^{\frac{l}{\lambda }}\\
&= \sum_{n=0}^\infty \left( \frac{1}{k!} \sum_{l=0}^k  {k \choose l} (-1)^{k-l} (l)_{n,\lambda } \right) \frac{t^n}{n!}.
\end{split}\end{equation}

Thus, by comparing the coefficients on both sides of \eqref{26}, we have
\begin{equation}\begin{split}\label{27}
\frac{1}{k!} \Delta^k (0)_{n,\lambda } = \frac{1}{k!} \sum_{l=0}^k {k \choose l} (-1)^{k-l} (l)_{n,\lambda } = \begin{cases}
S_{2,\lambda }(n,k)&\text{if}\,\,n \geq k,\\
0&\text{if}\,\,n<k.
\end{cases}
\end{split}\end{equation}

By \eqref{27}, we get
\begin{equation}\begin{split}\label{28}
\frac{1}{k!} \Delta^k (0)_{n,\lambda } = S_{2,\lambda }(n,k),\,\,\text{if} \,\,n \geq k.
\end{split}\end{equation}

From Theorem 7 and \eqref{28}, we obtain the following corollary.
\begin{cor}
For $n,k \in \mathbb{N}\cup \{0\}$ with $n \geq k$, we have
\begin{equation*}\begin{split}
B_{k,n}(x|\lambda ) = (x)_{k,\lambda }\sum_{m=k}^n {n \choose m} \beta_{n-m,\lambda }^{(k)}(1-x) \frac{1}{k!} \Delta^k (0)_{m,\lambda }.
\end{split}\end{equation*}
\end{cor}

Now, we observe that
\begin{equation}\begin{split}\label{29}
(1+\lambda t)^{\frac{x}{\lambda }}&= \Big( (1+\lambda t)^{\frac{1}{\lambda }}-1+1)^x = \sum_{k=0}^\infty {x \choose k} \Big( (1+\lambda t)^{\frac{1}{\lambda }}-1 \Big)^k \\
&= \sum_{k=0}^\infty (x)_k\frac{1}{k!} \Big( (1+\lambda t)^{\frac{1}{\lambda }}-1 \Big)^k \\
&= \sum_{k=0}^\infty (x)_k \sum_{n=k}^\infty S_{2,\lambda } (n,k) \frac{t^n}{n!}\\
&=\sum_{n=0}^\infty   \left( \sum_{k=0}^n  (x)_k S_{2,\lambda }(n,k) \right)
     \frac{t^n}{n!}.
\end{split}\end{equation}

On the other hand,
\begin{equation}\begin{split}\label{30}
(1+\lambda t)^{\frac{x}{\lambda }} = \sum_{n=0}^\infty { \frac{x}{\lambda } \choose n }
\lambda ^n t^n = \sum_{n=0}^\infty (x)_{n,\lambda }\frac{t^n}{n!}.
\end{split}\end{equation}

Therefore, by \eqref{29} and \eqref{30}, we obtain the following theorem.
\begin{thm}For $n \geq 0$, we have
\begin{equation*}\begin{split}
(x)_{n,\lambda } =  \sum_{k=0}^n (x)_k S_{2,\lambda }(n,k).
\end{split}\end{equation*}
\end{thm}

By Theorem 2.9, we easily get
\begin{equation}\begin{split}\label{31}
(x-k\lambda )_{i,\lambda } = \sum_{l=0}^i (x-k\lambda)_l S_{2,\lambda }(i,l).
\end{split}\end{equation}

From Theorem 2.6, we have the following theorem.
\begin{thm}
For $n,i \in \mathbb{N}$, with $i \leq n$, and $x \in [0,1]$, we have
\begin{equation*}\begin{split}
 \sum_{l=0}^i (x-k\lambda)_l S_{2,\lambda }(i,l) =\frac{1}{(x \oplus_\lambda  (1-x))^{n-i}} \sum_{k=i}^n \frac{{k \choose i}}{{n \choose i}} B_{k,n}(x|\lambda ).
\end{split}\end{equation*}
\end{thm}

\end{document}